\documentclass[letterpaper, 10 pt, conference]{ieeeconf}  

\IEEEoverridecommandlockouts       
\usepackage{algorithm}
\usepackage{algpseudocode}
\usepackage{graphicx}
\usepackage{float}
\usepackage{fancyhdr,bbm,amsmath,amssymb}
\usepackage{cite}
\newtheorem{problem}{Problem}

\title{\LARGE \bf
Optimal Routing of Modular Agents on a Graph 
}

\author{Karan Jagdale and Melkior Ornik
\thanks{K. Jagdale is with the Department of Aerospace Engineering and the
Coordinated Science Laboratory, University of Illinois Urbana-Champaign,
Urbana, IL 61801, USA. 
        {\tt\small karansj2@illinois.edu}}%
\thanks{M. Ornik is with the Department of Aerospace Engineering and the
Coordinated Science Laboratory, University of Illinois Urbana-Champaign,
Urbana, IL 61801, USA.
        {\tt\small mornik@illinois.edu}}%
}

\begin{document}

\maketitle
\thispagestyle{empty}
\pagestyle{empty}

\begin{abstract}
Motivated by an emerging framework of Autonomous Modular Vehicles, we consider the abstract problem of optimally routing two modules, i.e., vehicles that can attach to or detach from each other in motion on a graph. The modules' objective is to reach a preset set of nodes while incurring minimum resource costs. We assume that the resource cost incurred by an agent formed by joining two modules is the same as that of a single module. Such a cost formulation simplistically models the benefits of joining two modules, such as passenger redistribution between the modules, less traffic congestion, and higher fuel efficiency. To find an optimal plan, we propose a heuristic algorithm that uses the notion of graph centrality to determine when and where to join the modules. Additionally, we use the nearest neighbor approach to estimate the cost routing for joined or separated modules. Based on this estimated cost, the algorithm determines the subsequent nodes for both modules. The proposed algorithm is polynomial time: the worst-case number of calculations scale as the eighth power of the number of the total nodes in the graph. To validate its benefits, we simulate the proposed algorithm on a large number of pseudo-random graphs, motivated by real transportation scenario where it performs better than the most relevant benchmark, an adapted nearest neighbor algorithm for two separate agents, more than 85 percent of the time.
\end{abstract}

\section{INTRODUCTION}
Modular systems --- systems where agents can attach and detach from others mid-mission --- are an emergent generalization of multi-agent systems. Their abstract formulation is motivated by the novel technology of Autonomous Modular Vehicles (AMV), which allows two vehicles to attach and detach from each other \cite{enwiki:1129248893}. 

Applications of modular agents are primarily studied for public transportation to lower the operation cost of vehicles and improve the service quality to customers. For example, modular bus systems are proven to show significant benefits over the non-modular bus system \cite{khan2023application}. Upon joining, these buses have an open area where the passengers can stand and travel from one bus to another. This ability to distribute the passengers between two buses (which can detach later) decreases the average passenger travel time. Vehicle modularity has been studied in specific transportation scenarios like oversaturated traffic --- where the passenger demand is higher than the transportation network capacity \cite{chen2020operational}, flex route transit service \cite{liu2021improving}, shared-use corridors where different bus routes sharing a common bus stop \cite{shi2020variable} etc. 
 Papers \cite{pei2021vehicle, ji2021scheduling} utilize the modularity to vary the vehicle capacity and propose an optimal vehicle dispatch schedule given the passenger demand. 

In addition to the public transportation domain, modular vehicles can play a crucial role in improving vehicle platooning --- a method of driving a group of vehicles together to increase the road capacity --- for enhanced mobility and passenger comfort, as studied by paper \cite{li2022trajectory}. Vehicular modularity has also been considered in military vehicles: by attaching two military vehicles of different types, a unit can obtain increased mobility on a variety of terrains \cite{dasch2016survey}.
\subsection{Motivating Example} \label{sub:motex}
\begin{figure}[H]
    \centering
    \includegraphics[scale = 0.25, angle = 90]{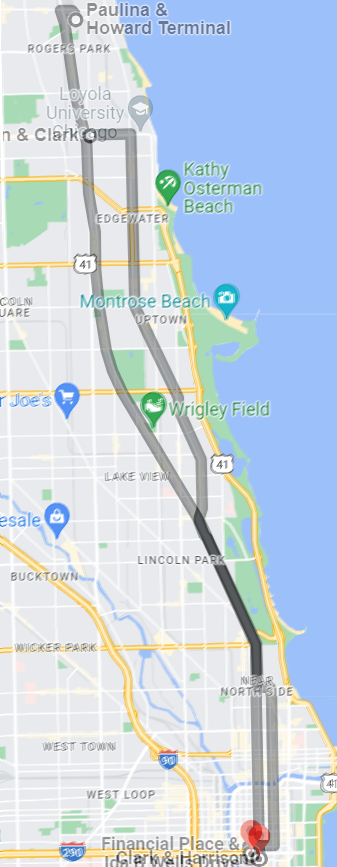}
    \caption{Map showing routes of CTA Bus $22$ South and $36$ South \cite{cta} in grey highlight. Starting locations of the buses are towards the left border of the figure. The dark grey highlight depicts the route where the buses share common stops.}
    \label{fig:cta}
\end{figure}

Vehicle modularity can be used to improve the existing bus transportation system. To exemplify, we present bus routes of two bus lines in Chicago, namely $22$ South and $36$ South. These buses depart from locations relatively close to each other, and their routes are initially side by side. Their paths merge for some part and eventually get separated, as shown in Fig. \ref{fig:cta}. The buses have common stops on the part where their routes are merged. Modularity can prove useful in this case by joining the buses in the region where their paths are merged. Upon joining the buses, the passengers in the joined bus have access to the bus stops of both buses, without the need to change buses. As the passengers decide which of the two buses to sit on, they save waiting time at the bus stop. Additionally, the passengers boarding the joined bus will also have access to the bus stops of both buses as they can choose which of the two buses (of the joined bus) to enter. In the joined state, the buses will have higher fuel economy and contribute less to traffic congestion.
\subsection{Related work and contributions}
Previous efforts in modular agents are primarily focused on the passenger transport system and are developed for specific problems where the agent routes are fixed. To the best of our knowledge, modular vehicles are not studied in an abstract setting. This paper considers the problem of routing modular agents on a graph to reach preset nodes (targets) by traversing the minimum possible distance. In the rest of the discussion, the module refers to the agent that cannot split further, and the agent is any vehicle possibly formed by joining multiple modules. We consider the scenario where the cost of traversing an edge for an agent is the same as that of a single module to simplistically model benefit of joining two modules like fuel efficiency \cite{li2022trajectory}, less traffic congestion (due to less number of vehicles). 

This paper proposes a novel algorithm that routes two modules to the preset targets. While the problem of doing so with minimal cost is NP-hard, we present a heuristic approach which requires polynomial time in the number of targets and graph size.  The algorithm aims to determine when and where to join or split the modules by predicting the cost incurred to complete the mission by a single module vs. two modules starting at multiple nodes. Intuitively, the modules are joined if they are relatively close to each other and away from the remaining targets, and the joined agent splits once it is routed enough close to the remaining targets. 

The paper's outline is as follows:  We state the problem formally and present an example to demonstrate the potential of modularity in Section \ref{sec:ps}. In  Section \ref{sec:tnn}, we present the adapted nearest neighbor algorithm for two agents which is used to formulate the proposed algorithm and also used as a benchmark for comparison of performance of the modular agents with non-modular agents. We propose an optimal routing algorithm in the Section \ref{sec:optalg} and analyze its complexity in Section \ref{sec:compan}. We present a class of graphs for which the proposed algorithm \textit{provably} performs better than non-modular agents in Section \ref{sec:poo}. To illustrate the algorithm, we present results obtained using the proposed algorithm on a large number of pseudo random graphs, inspired by real transportation scenario in Section \ref{sec:res}.

\section{PROBLEM STATEMENT} \label{sec:ps}

This paper considers the problem of optimal planning with modular agents traversing on an undirected graph $\mathcal{G}(\mathcal{V}, \mathcal{E})$. The modules aim to visit a set of  preset nodes while incurring the least cost of travel. These preset nodes are referred to as \textit{targets} in the paper, and the set of the targets is given by $\mathcal{T} \subset \mathcal{V}$. The cost incurred by a module on its path is equal to the sum of the weights of the traversed edges, while the agent formed by joining two modules incurs the same cost as the individual modules. Exploiting their modular capabilities, modules can join and split at any node in the graph. Moreover, a module can also choose not to move at a given time.
\par Let us mathematically define the above problem of optimal planning for modular agents. The sequence of nodes traversed by a module defines its \textit{path}. For every edge $e\in\mathcal{E}$, let its weight, i.e., the cost incurred by traversing it, be denoted by $w_e>0$. Let the set of modules moving at time instant $t$ is given by $\mathcal{K}(t)$ and let the edge traversed by Module i at that time be instant denoted by $e_i(t)$, then the total cost incurred by modules during the mission is given by
\begin{align}\label{eqn:tc}
    \sum_{t=1}^T\sum_{e \in \cup_{i\in \mathcal{K}(tz)}\{e_i(t)\}}w_e
\end{align}
where T is the time step after which the modules' mission is complete. In equation \ref{eqn:tc}, we assume that all modules that traverse the same edge at the same time are joined. Given the proposed cost model, there is no benefit in modules not joining when traversing the same edge concurrently. The actions of joining and splitting are thus automatically encoded by keeping a record of timed paths of individual modules. 
\begin{problem} \label{p1}
    Let $n$ modules operate on a graph $\mathcal{G}(\mathcal{V,E})$ with a target set $\mathcal{T}\subset \mathcal{V}$. Denote the path of Module $i$ by $P_i=(v_i(0),...,v_i(T))$  with $v_i(t)\in\mathcal{V}$ and $(v_i(t),v_i(t+1))\in\mathcal{E}$ for each $0\leq t<T$. Determine paths $P_1, ..., P_n$ which minimize \ref{eqn:tc} such that $\mathcal{T}\subset\cup_{i=1}^n\cup_{t=0}^T\{v_i(t)\}$
\end{problem}

\begin{figure}[H]
    \centering
    \includegraphics[scale = 0.4]{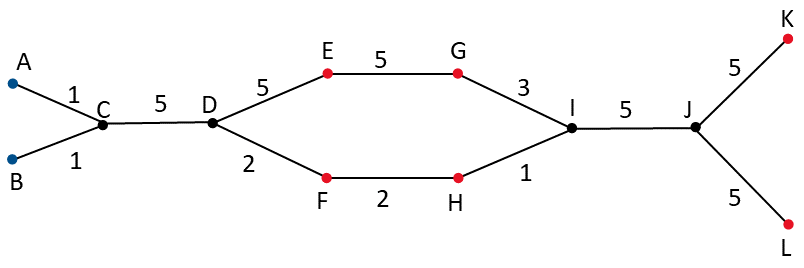}
    \caption{Example graph with modules' initial locations shown with blue dots and the targets indicated with the red dots. The numbers next to the edges indicate the costs of traversing the edges.}
    \label{fig:mot_schm}
\end{figure}

Problem \ref{p1} is in general NP-hard, as the case of n=1 reduces to the classical NP-complete problem of finding a minimum Hamiltonian path \cite{garey1974some}. As a first step in finding a computationally feasible approximate solution, this paper considers the case of $n=2$. Such a scenario already exhibits the fundamental characteristics of modularity by allowing the two modules to join into a combined agent and subsequently split. To demonstrate the potential benefits of modularity in such a case, we provide a short example illustrated in Fig. \ref{fig:mot_schm}. The two modules start from nodes $A$ and $B$, and need to visit target set $\mathcal{T}=\{E,F,I,J\}$. The graph structure and weights  are given on Fig. \ref{fig:mot_schm}. The optimal paths for two modular agents are given by $(A,C,D,E,G,I,J,K)$ and $(B,C,D,F,H,I,J,L)$. The modules join at node $C$, split at $D$, join again at $G$, and finally split at $H$. The total cost incurred is $40$. If the agents were not modular, i.e., lacked the capability to join and split, by inspection we can verify that the optimal routing would be that one of the modules does not move at all, whereas the other visits all the targets, e.g., $(A,C,D,F,D,E,G,I,H,I,J,K,J,L)$. The cost incurred in this case is $45$.

We now move to propose an approximately optimal solution to Problem 1 in the case when $n=2$.

\section{Routing algorithm} \label{sec:optalg}
 As described earlier, the analytical solution to Problem \ref{p1} is NP-hard; we opt for heuristic approach. Our solution consists of three interacting elements which are considered anew at every time step:
\begin{enumerate}
    \item If the modules are separated, deciding if and where to join them: We do that by comparing distance between the modules with the distance of the remaining targets from the modules. Intuitively, we join the modules if they are close enough to each other and away from the targets.
    \item If the modules are joined, deciding if and where to split them. We make the decision this by estimating the costs of paths obtained by splitting the joined agent at different nodes, including not splitting at all, and visiting the remaining targets using the adapted nearest neighbor algorithm for two agents. We find the optimal splitting node that corresponds to the minimum cost.
    \item Optimal routing of the modules towards the targets. We decide whether to route one of the modules or both of them at the current time step by predicting the cost of doing so using nearest neighbor strategy.
\end{enumerate}
 Now, we present a detailed analysis of the above decisions and the tools used within one by one; we start with the first ingredient in our computationally tractable solution: an adapted nearest neighbor algorithm for two concurrently moving agents.
\subsection{Adapted nearest neighbor for two agents} \label{sec:tnn}
As we described above, the problem of modular optimal planning is a generalization of the minimum Hamiltonian path problem. A nearest neighbor approach \cite{nilsson2003heuristics}, in which an agent always visits the closest unvisited target, has been widely used to approximate a solution to this problem for a single agent. We now adapt this approach for two agents, and use it in the formulation of the proposed algorithm and later as a benchmark of optimal policies of \textit{non-modular agents} for comparison with the proposed algorithm.

The adapted nearest neighbor algorithm for a target set $\mathcal{T}$ and module nodes $a_1$ and $a_2$ is as follows: Let $d(n_1,n_2)$ be the minimal possible cost of a path between nodes $n_1,n_2$. First, we find $(\tau_1, \tau_2)\in \mathcal{T}^2$ with $\tau_1\neq \tau_2$ such that $d(a_1,\tau_1)+d(a_2, \tau_2)=\min_{\tau_1'\neq \tau_2'} d(a_1,\tau_1')+d(a_2,\tau_2')$. Then we move Agent 1 and Agent 2 by one node on the shortest path to $\tau_1$ and $\tau_2$ respectively. While we consider $n=2$ in keeping with the remainder of the paper, this algorithm can be directly generalized for any $n>1$.


\subsection{Joining decision}
\begin{figure}[H]
    \centering   \includegraphics[scale = 0.35]{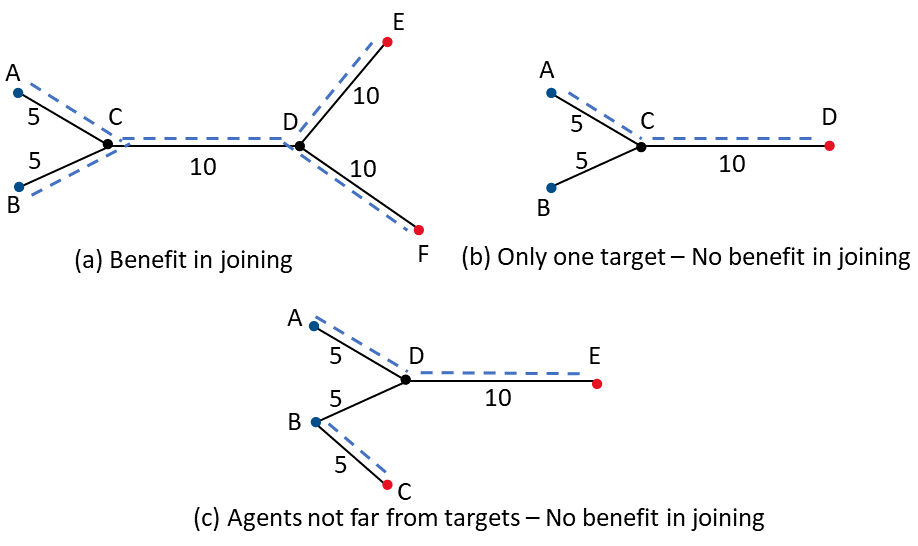}
    \caption{Different scenarios to explain the reasoning behind the benefits of joining. Blue dots show the modules' initial locations, and the red dots show the target locations. Blue dashed lines show the optimal module path.}
    \label{fig:jndec}
\end{figure}

Intuitively, joining the modules is beneficial if the modules are relatively close to each other and far from the remaining targets and there are sufficiently many targets remaining to be visited. This intuition is illustrated by Fig. \ref{fig:jndec}. In case (a), joining the two modules is beneficial, and the cost of the optimal path, as indicated in the above figure, is $40$. If the modules did not join, only one of the modules would visit both targets, or both modules would visit one of the targets, and the cost incurred in these cases would be $45$ and $50$, respectively. In case (b), joining the modules is not beneficial as only one target is remaining and there is no need for both modules to move. In case (c), module locations are not sufficiently away from the target nodes. 
\par To formalize the notion of the modules being far from the remaining targets, we will define a set of \textit{central nodes} $\mathcal{C}$ with respect to the nodes where the modules are present. $\mathcal{C}$ is given by
\begin{equation*}
    \mathcal{C}(t) = \arg\min_{v \in \mathcal{V}} \frac{d(v,m_1(t)) + d(v,m_2(t))}{2},
\end{equation*}
where $m_1(t)$ and $m_2(t)$ denotes Module 1 and Module 2 nodes at time $t$. Let $d_c(t)$ be the distance between central node $c \in \mathcal{C}(t)$ and the closest remaining target, i.e.,

\begin{align}
     d_c(t) = \min_{\tau \in \overline{\tau}(t)} d(c, \tau),      \label{eqn:agTar}
\end{align}
where $\overline{\tau}(t)$ denotes the set of targets remaining at time $t$.
We formulate the following joining condition: if $d(m_1(t),m_2(t))<d_c(t)$, then direct the modules to join at node $n_j$ given by 
\begin{equation}
    n_{j} = \arg\min_{c \in \mathcal{C}} d_c(t). \label{eqn:jnnd}
\end{equation}
In the case $n_j$ is not uniquely defined, we choose any node in the set defined by \ref{eqn:jnnd}. The above formulation ensures that we join the modules at the node which is both central with respect to the modules' positions and ``in the direction of" the remaining targets. Thus, the modules move closer to the remaining targets while moving towards the joining node.

\subsection{Splitting decision}\label{sub:split}
\begin{figure}[H]
    \centering
    \includegraphics[scale = 0.35]{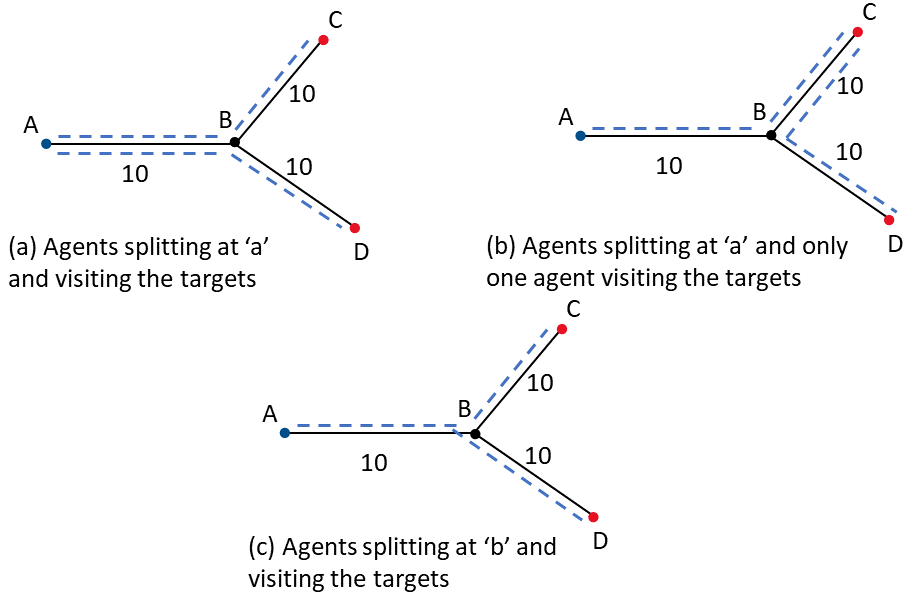}
    \caption{Different scenarios depicting intuitive reasoning behind the optimal split location. Blue dots show the modules’ initial locations, and the red dots show the target locations. Blue dashed lines show the modules'/agent's path.}
    \label{fig:spnd}
\end{figure}

To formulate a splitting condition, we use a similar intuition to that of joining.
Routing the joined agent close to the remaining targets and splitting there will produce a lower cost than splitting first and routing the modules separately. To exemplify this, we refer to Fig. \ref{fig:spnd}. In strategy (a), modules split at node A, and both modules visit the targets. The cost of routing obtained in this case is $40$. In strategy (b), the agent splits at node A, and only one of the modules visits the remaining targets. The cost of routing is again $40$. However, in strategy (c), the joined agent first moves close to the remaining targets, i.e., moves to node B and splits there, and each module visits a target. The routing cost obtained in this case is $30$. 

With the above intuition, we employ the following method to determine the optimal splitting node.
\par For each $r \geq 0$ We calculate the cost of the agent traveling for $r$ nodes on its nearest neighbor path --- path obtained using the nearest neighbor approach for the agent --- followed by splitting the agent into two modules, with the separated modules  proceeding using the adapted nearest neighbor algorithm for two agents. 
We choose the splitting node corresponding to the $r$, which produces smallest cost. 
 
\par In a case where the targets are naturally clustered, the above strategy might not produce optimal results because the modules will split while in the last cluster to be visited, and not at any earlier point which is not optimal as demonstrated in Section \ref{sec:ps}. To address such cases, we include an additional clustering filter to the splitting decision. We divide the targets into two clusters using the standard technique of $k$-means clustering \cite{hartigan1979algorithm}; if the \textit{DB-index} \cite{davies1979cluster} of the clusters is lover than a threshold value, i.e., the clusters are fairly separated, we only consider the targets in the cluster closest to the agent location to obtain optimal splitting node. 

\par While the intuitions for the joining condition and for splitting are the same, the technical approach differs. Using analogy of joining, one can get optimal splitting node as the one which is central with respect to the targets and closest from the agent location. Although, for the case when the number of targets is greater than two, this method will not produce optimal splitting node. Consider a case with four target nodes arranged to look like vertices of a quadrilateral, the most central node for the target nodes will be at the node closest to the centroid of the quadrilateral. Now, routing the agent all the way to the most central node and splitting there in the above case is clearly not optimal.

To avoid the possible infinite loop of performing repeated joining and splitting, we do not let the agent join immediately after it is split. Instead, in the next step, we route the agent as given in Section \ref{sub:taras}.

\par The splitting strategy is summarized in \textbf{Algorithm} \ref{alg:optsp}. Note that, in the pseudocodes described below, $p(n_1, n_2)$ gives the sequence of nodes corresponding to the shortest path (least cost) between $n_1, n_2$ starting with $n_1$ and $d(n_1, n_2)$ gives the shortest distance(least cost) between nodes $n_1, n_2$.
\begin{algorithm}[H]\label{a:alg1}
\caption{Node to split the agent}\label{alg:optsn}
\label{alg:optsp}
 \hspace*{\algorithmicindent} \textbf{Input:} Agent node : $a$, target set: $\mathcal{T}$\\
 \hspace*{\algorithmicindent}        Threshold \textit{DB-index}: $D_t$
 \\
 \hspace*{\algorithmicindent} \textbf{Output:} Node to split the agent $n_{sp}$
\begin{algorithmic}[1]

\State Divide $\mathcal{T}$ into two clusters $\mathcal{C}_1, \mathcal{C}_2$ with $\mathcal{C}_1$ closer to agent.
\State $D$: \textit{DB-index} of $\mathcal{C}_1, \mathcal{C}_2$ 
\If{$D<D_t$}
    \State $\overline{\tau} = \mathcal{C}_m$
\Else
    \State $\overline{\tau} = \mathcal{T}$
\EndIf
\State $P$ : Agent's nearest neighbor path to visit targets in $\overline{\tau}$
\For{r=0:length(P)-1}
    \State $n_{r} = P(r)$ 
    \State Compute $cost_r$, i.e., cost  of visiting remaining   
    \State nodes using adapted nearest neighbor with modules 
    \State starting from $n_r$
    \State $c_r = d(a,n_r) + cost_r$
\EndFor
\State $p = \arg\min_r(c_r)$
\State For any $\overline{p} \in p, n_{sp} = P(\overline{p})$
\end{algorithmic}
\end{algorithm}
Now we will analyze the decision of routing one Vs. both the modules.

\subsection{Target assignment}\label{sub:taras}
To understand the significance of this section, consider the scenario in Fig. \ref{fig:mot_schm} with targets at only F and J nodes. In this case, it is not optimal to join the modules as there is no benefit in splitting later as the shortest path to the farther target passes through the closer target. The optimal solution would be only one of the modules going to both targets. 

To decide the number of modules to be routed, we compare the cost incurred by each module if it visits all the targets using the nearest neighbor approach with the cost of visiting all the targets using both modules with the adapted nearest neighbor for two agents. Based on the computed costs, we move the modules by one node using the strategy that produces the least cost and compute the optimal route. Combining the joining, splitting, and target assignment mechanisms, we obtain the complete routing algorithm, presented in \textbf{Algorithm} \ref{alg:mainalg}.

\begin{algorithm}[H]\label{a:alg2}
\caption{Optimal Routing}
\label{alg:mainalg}
 \hspace*{\algorithmicindent} \textbf{Input:} Graph with the module nodes : $m_1, m_2$ and \\
 \hspace*{\algorithmicindent} the target nodes: $\mathcal{T}$ \\
 \hspace*{\algorithmicindent} \textbf{Output:} Modules' timed paths : $P_1, P_2$
\begin{algorithmic}[1]
\State Initialize $P_1(0) = m_1, P_2(0) = m_2, t=0$
\State lastSplit = False
\While{All targets are not visited}
\If {$P_1(t) = P_2(t)$}
\State Obtain splitting node $n_{sp}$ using \textbf{Algorithm} \ref{alg:optsn} 
\State $m_1 = m_2 = n_{sp}$
\State $P_1 = [P_1, p(m_1,n_{sp})], P_2 = [P_2, p(m_2, n_{sp})]$
\State $t = length(P_1)$
\State lastSplit = True
\Else
    \If{$ d(m_1, m_2) \leq d_c$ and lastSplit = False}
        \State $P_1 = [P_1,p(m_1,n_{jn})]$,
        \State $P_2 = [P_2, p(m_2, n_{jn})]$
        \State $n_1 = length(P_1), n_2 = lenght(P_2)$
        \If{$n_1 < n_2$}
            \State $P_1 = [P_1, n_{sp}\times ones(1,n_2-n_1)]$
            
        \Else
            \State $P_2 = [P_2, n_{sp}\times ones(1,n_1-n_2)]$
        \EndIf
        \State $t = length(P_1)$
        
    \Else
        \State Compute:
        \State $p_1$ : Cost of visiting all the targets with  
        \State Module 1 using nearest neighbor algorithm
        \State $p_1$ : Cost of visiting all the targets with  
        \State Module 2 using nearest neighbor algorithm
        \State $p_3$ : Cost of visiting all the targets with   
        \State both modules using adapted nearest neighbor.
        \State $i_m = \arg\min_i(p_i)$
        \If {$i_m = 1$ or $i_m = 2$}
            \State $\tau^c_{i_m} = \arg \min_{\tau \in \mathcal{T}}d(m_{i_m}, \tau)$
            \State $P = p(m_{i_m}, \tau^c_{i_m})$
            \State $m_{i_m} = P(2)$
            \State $P_{1} = [P_{1}, m_1], P_2 = [P_2, m_2] $
        \Else
            \State Obtain module paths $N_1, N_2$ using 
            \State adapted nearest neighbor.
            \State $m_1 = N_1(2), m_2 = N_2(2)$
            \State $P_1 = [P_1, m_1], P_2 = [P_2, m_2]$
        \EndIf       
        \State $t= t+1$
    \EndIf
    \State lastSplit = False
\EndIf
\EndWhile
\end{algorithmic}
\end{algorithm}

\section{Complexity analysis} \label{sec:compan}
We will show that Algorithm 2 for modular planning both generally outperforms the benchmark of planning for two non-modular agents --- provably so for a particular class of graphs --- and does it in a computationally feasible manner. We begin by analyzing its computational complexity. We separately investigate the cases of modules being separated and joined; then, we state the overall complexity of the algorithm. 
\par Let $m,n$ denote the total number of nodes and the total number of targets, respectively. We use Dijkstra's algorithm to find the shortest path between two nodes, which has the time complexity of $O(m^2)$ \cite{barbehenn1995efficient}. In the nearest neighbor algorithm for a single module, at each step, we find the closest target by computing the shortest distance from the module to $n$ targets using Dijkstra's algorithm. Thus, the complexity of finding the nearest target is $O(n m^2)$. Until the module visits all the targets, we repeat the process of finding the closest target $n$ times. Thus, the complexity of using the nearest neighbor to visit the targets with a single module is $O(n^2m^2)$. Similarly, the complexity of the adapted nearest neighbor algorithm is also $O(n^2 m^2)$, as it finds the target nearest to the modules $4$ times performing $O(n m^2)$ computations, and this process is repeated for $n/2$ times to visit all the targets. 

The complexity of deciding whether to join or not is $O(nm^2+m^3) = O(m^3)$. Namely, we first compute the central node with respect to the module locations, which requires computing the shortest distance from each node to the module locations, i.e., total $O(2m \times m^2) = O(m^3)$ calculations. Then, we obtain the closest target from the central node by performing $O(n \times m^2) = O(nm^2)$ calculations. We perform computations to decide whether to join or not at every time step when the modules are separated.

When the modules are separated and we did not decide to join them in future, we move them by one node. To decide the next nodes to move the modules, we compute the nearest neighbor path individually for both modules and the path using adapted nearest neighbor. As described earlier, the complexity of these operations is $O(n^2 m^2)$. Thus, the complexity of deciding the next step when modules are separated is $O(n^2m^2 +  m^3)$. 
\par When the modules are joined, they first traverse using the nearest neighbor strategy by $r$ nodes and then visit the remaining targets if any separately using the adapted nearest neighbor. As described in Section \ref{sub:split}  optimal splitting node depends on the value of $r$ which produces the least cost. The number of computations required to find the agent's nearest neighbor path is $O(n^2m^2)$ and the number of computations for doing adapted nearest neighbor from the $r^{th}$ node is also upper-bounded by $O(n^2m^2)$. The number of times adapted nearest neighbor 
 is implemented to compute the optimal splitting node is upper-bounded by $m$. Thus, the computations required to find the optimal splitting node are upper-bounded by $O(m^2n^2 + n^2m^3)$, i.e., the complexity of finding the splitting node is $O(n^2m^3)$.
 
 Thus,  at each time step, we perform $O(n^2m^3)$ computations if the modules are joined, whereas, we perform  perform $O(n^2m^2 + m^3)$ computations if the modules are separated. The number of time steps to visit all the targets is upper-bounded by $O(m^2n)$: the same pair of nodes cannot be visited more than twice without visiting a target in-between. Thus the worst case complexity of the proposed algorithm is $m^2n \times O(n^2m^3) = O(n^3m^5)$, i.e., given $n \leq m, O(m^8)$.
\section{Performance on Clustered Graphs}\label{sec:poo}
As a first step for a future theoretical discussion of our proposed algorithm's performance, we present the class of graphs where the proposed algorithm with modular agents will \textit{provably} always produce lower cost than the non-modular agents. 

\begin{figure}[H]
     \centering
     \includegraphics[scale = 0.4]{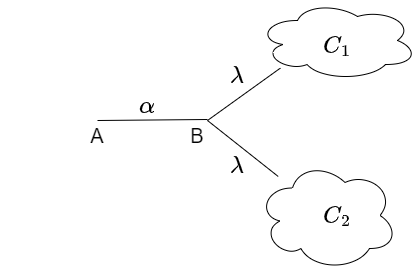}
     \caption{Figure describing the class of graphs where the proposed algorithm with modular agents are guaranteed to produce lower cost than non-modular agents.}
     \label{fig:graphclass}
 \end{figure}

We consider a class of graphs where the targets can be categorized into clusters. First, we define the notion of cluster; $\mathcal{C} \subset \mathcal{V}$ is a cluster is it satisfies following condition:
\begin{align*}
    \max_{c,c' \in \mathcal{C}}d(c,c') < \min_{c\in \mathcal{C}, b \notin \mathcal{C}} d(c,b).
\end{align*}
Fig. \ref{fig:graphclass} illustrates the considered class of graphs. The clusters are shown by $\mathcal{C}_1$ and $\mathcal{C}_2$. While this class is obviously simplistic --- e.g., it has only two clusters, , only two nodes that are not in the clusters, and the distances of the clusters to those nodes is equal --- the proof below can be directly extended to graphs with slightly more complicated graphs. Formally defining such graphs would require burdensome notation, so we omit such a discussion. While noting that the same proof intuition continues to hold. On the other hand, extending this proof further to \textit{substantially more complicated graphs} is one of the central objectives for future work.
  
Modules in the considered class of graphs are initially joined and are at node \textbf{A}. Nodes \textbf{A} and \textbf{B} are $\alpha$ distance apart. Clusters $C_1$ and $C_2$ are $\lambda$ distance away from node \textbf{B}, i.e., 
\begin{align*}
    \lambda = \min_{\tau \in  \mathcal{C}_1} d(B, \tau) = \min_{\tau \in  \mathcal{C}_2} d(B, \tau) 
\end{align*}
Let $\beta_1$ and $\beta_2$ be the cost of visiting all targets in $C_1$ and $C_2$ starting from the node closest to \textbf{B} in $\mathcal{C}_1$ and $\mathcal{C}_2$, respectively, using the nearest neighbor approach.

Let us compare the performance of modular agents to a non-modular benchmark. First, consider the case of non-modular agents: Depending on the relationship between alpha and lambda, one of the two following strategies is optimal:
 \begin{enumerate}
     \item Either Agent 1 or Agent 2 alone visits all the targets, and the cost incurred is given by
     \begin{align}\label{eqn:nm1}
          \alpha + 3 \lambda + \beta_1 + \beta_2.
     \end{align}
     \item Each agent visits the targets in one of the clusters, and the cost incurred is given by,
     \begin{align}\label{eqn:nm2}
          2 \alpha + 2 \lambda + \beta_1 + \beta_2.
     \end{align}
 \end{enumerate}
Now, let us consider the proposed algorithm for modular agents. The algorithm will predict the cost of visiting all the targets by splitting at different nodes on its nearest neighbor path i.e., it will compute the cost of splitting at
\begin{enumerate}
    \item \textbf{A}: The cost will be equal to \ref{eqn:nm2} as the modules use the adaptive nearest neighbor after splitting.
    \item \textbf{B}: $ \alpha + 2 \lambda + \beta_1 + \beta_2$, i.e, each module visiting the targets in one of the clusters, as the clusters are equidistant from \textbf{B} and consequently the adaptive nearest neighbor assigns one target in each cluster to both modules when the agent is at node B.
    \item any node in $\mathcal{C}_1$ or $\mathcal{C}_2$: The cost will be roughly equal to $\alpha + 3 \lambda + \beta_1 + \beta_2$, as both modules first visit targets in the clusters they entered and then visit targets in the other cluster. The cost of routing two modules with the adaptive nearest neighbor is slightly different than the cost of routing using the nearest neighbor using a single module in the clusters.  
\end{enumerate}
As the predicted cost of splitting at \textbf{B} is lowest, the agent will be routed to node \textbf{B}, paying the cost of alpha, and split there. 
The algorithm will then evaluate the cost of visiting all the targets with
\begin{enumerate}
    \item only Module 1 or Module 2, i.e., $3 \lambda + \beta_1 + \beta_2$,
    \item using both modules with the adaptive nearest neighbor, i.e., $2 \lambda + \beta_1 + \beta_2$.
\end{enumerate}
Once the modules are inside the clusters, they will be routed using adaptive nearest neighbor, which basically is the nearest neighbor on each module separately as the modules are in different clusters . Thus, the total cost using the proposed algorithm is given by
 \begin{align}\label{eqn:m}
     \alpha + 2 \lambda + \beta_1 + \beta_2.
 \end{align} 
 Comparing \ref{eqn:nm1}, \ref{eqn:nm2} with \ref{eqn:m} makes it clear that the cost obtained with the proposed algorithm using the modular agents is always less than the cost obtained with the non-modular agents.
\section{Numerical Results}\label{sec:res}
To describe the performance of the proposed algorithm, we consider an illustrative example motivated by real-life scenario in in retail logistics industry and simulation on a batch of random graphs. 

\subsection{Illustrative Example}
\begin{figure}[H]
    \centering
    \includegraphics[scale = 0.5]{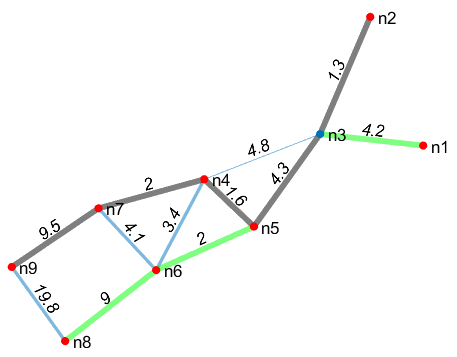}
    \caption{Graph representing the road network connecting $8$ warehouses. The number next to the edges is the distance in miles between the locations represented by the nodes. Green and grey highlights indicate the path of Truck 1 and Truck~2, respectively. Red dots highlight the warehouse locations.}
    \label{fig:spjn4}
\end{figure}

\begin{figure}[H]
    \centering
    \includegraphics[scale = 0.5]{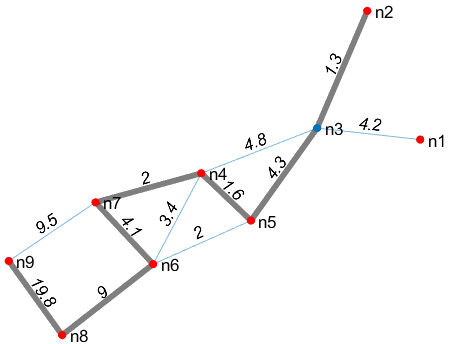}
    \caption{Graph showing Trucks' path in green and grey highlight. Note that, Truck~1 corresponding to green highlight does not move. Warehouse locations are highlighted with red dots. The number next to the edges is the distance in miles between the locations represented by the nodes.}
    \label{fig:nspjn4}
\end{figure}

in retail logistics, companies use warehouses to store manufactured goods before distributing them to shops. Usually, big companies have a network of warehouses spread near big cities.
We consider the problem of distributing shipments in a set of Amazon warehouses located in/near Los Angeles, California, USA. We create a graph using the data from google maps \cite{gmap} from several warehouse locations and use the proposed algorithm to obtain the optimal path using two modular trucks. The trucks start from two warehouses and want to visit another $6$ warehouses to deliver the shipments. Fig. \ref{fig:spjn4} shows the graph representing the scenario along with the optimal truck routes.
 As shown in Fig. \ref{fig:spjn4}, the trucks join at \textbf{n3} and split at node \textbf{n5}. The total cost incurred is $33.9$. On the other hand, Fig. \ref{fig:nspjn4} shows the optimal routes when they are non-modular, i.e., two separate trucks without the capability to join or split. In this case, only the truck starting at \textbf{n2} moves to visit the warehouses, and the cost incurred is $42.1$.
\subsection{Batch simulation on random graphs} 
\begin{figure}[H]
    \centering
    \includegraphics[scale = 0.5]{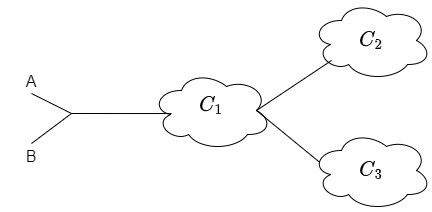}
    \caption{Figure indicating the considered class of graphs where the target nodes are present in three clusters. The modules' initial locations are marked by A, B.}
    \label{fig:simclus}
\end{figure}

Moving to show the comparative performance of the proposed algorithm on a large set of graphs, we run $100$ simulations on graphs generated in a pseudo-random fashion having a total of $18$ nodes and $8$ target nodes. These pseudo-random graphs are produced in a way that ensures that target nodes are located in three clusters, with a structure shown in Fig. \ref{fig:simclus}. Such a structure is motivated by the public transportation example described in Section \ref{sub:motex}, where multiple buses partially go through the same area. We use $0.4$ as the threshold value of \textit{DB-Index} ($D_t$) in the simulations. Note that the considered class of graphs is more general than those in Section \ref{sec:poo}, for which we provided a theoretical proof of performance. 
Intuitively, having a cluster structure ensures that the split and join actions can reduce the total routing cost. Although we cannot guarantee the performance of the proposed algorithm on the graph structure shown in Fig. \ref{fig:simclus}, the proposed algorithm indeed produced a smaller cost than routing for non-modular agents in $85$ out of $100$ graphs.

\section{CONCLUSIONS} \label{sec:conc}
Spurred by the emergent framework of Autonomous Modular Vehicles, this paper presents the first known formulation of an optimal routing problem for modular agents on a graph. In particular, the paper considers the problem of optimal planning for modules to visit a set of preset nodes while incurring minimum cost on a graph. We assume that the cost incurred by a joined agent --- a vehicle formed by joining two modules --- is the same as that of the single module. This cost formulation implicitly captures the benefit of joining two vehicles, like less traffic congestion and higher fuel economy. 

We propose a heuristic approach to solve the problem as the analytical solution to the problem is NP-hard. 
We present a routing algorithm for two modules as a starting step. Our formulation relies on three decision making elements: Joining, Splitting and routing. We join the modules if they are close to each other and away from the remaining targets and enough targets remain to be visited. For splitting, we compute the cost of visiting all the targets for the cases where the agent splits every node on its nearest neighbor path one by one, and the modules complete the routing using the nearest neighbor algorithm adapted for two agents. We choose the node to split the agent corresponding to the minimum routing cost. When the modules are neither joining nor splitting, we direct one or both of them towards remaining targets, depending on which option incurs less cost, using a method based on the nearest neighbor approach: visiting all the targets using the nearest neighbor with only one of the modules, visiting the remaining targets using adapted nearest neighbor for two agents.

The motivating application domain of modular-agent public transportation drives challenges remaining for future work. Namely, a realistic public transportation scenario will consist of multiple agents and modules, with agents possibly composed of more than two modules, operating on a large and intricately structured graph. Developing such a framework presents both a challenge in terms of algorithm design and computational cost because increasing the number of modules increases the number of decisions to make at each step. The proposed algorithm with modular agents is proved to be superior than non-modular agents for a simplistic class of graphs in the paper. Extending this proof to substantially complicated graphs is one of the primary goals of the future work. Additionally, the current resource cost model, in which the resource cost of the joined agent is the same as that of the separate module, is simplistic. The application domain consists of multiple objective functions, including resource consumption, contribution to traffic congestion, and travel time. One can include the traffic condition in the cost formulation to optimize the system performance for the above objectives. Finally, the proposed algorithm should be validated in a high-fidelity traffic simulator and improved, keeping in mind the likely stochastic nature of passenger arrival and travel time of the modules/agents.

\addtolength{\textheight}{-12cm}   




\bibliographystyle{plain} 
\bibliography{ref}

\end{document}